\documentclass{amsart}
\usepackage{amssymb,amsmath}
\newtheorem{Theorem}{Theorem}[section]
\newtheorem{Definition}[Theorem]{Definition}
\newtheorem{Proposition}[Theorem]{Proposition}

\newtheorem{Lemma}[Theorem]{Lemma}
\newtheorem{Corollary}[Theorem]{Corollary}
\theoremstyle{remark}
\newtheorem{Example}[Theorem]{Example}

\def\eps{\varepsilon}

\def\supp{\operatorname{supp}}

\def\be{\begin{enumerate}}
\def\ee{\end{enumerate}}
\def\bT{\begin{Theorem}}
\def\eT{\end{Theorem}}
\def\bP{\begin{Proposition}}
\def\eP{\end{Proposition}}
\def\bD{\begin{Definition}}
\def\eD{\end{Definition}}
\def\bE{\begin{Example}}
\def\eE{\end{Example}}
\def\bL{\begin{Lemma}}
\def\eL{\end{Lemma}}
\def\bC{\begin{Corollary}}
\def\eC{\end{Corollary}}

\def\J{{\mathcal J}}

\begin{document}
\title{Green functions, the fine topology and restoring coverings}
\author{Tony L. Perkins}

\address{ Department of Mathematics,  215 Carnegie Building,
Syracuse University,  Syracuse, NY 13244-1150}
\begin{abstract}
There are several equivalent ways to define continuous harmonic functions $H(K)$ on a compact set $K$ in $\mathbb R^n$. One may let $H(K)$ be the unform closures of  all functions in $C(K)$ which are restrictions of harmonic functions on a neighborhood of $K$, or take $H(K)$  as the subspace of $C(K)$ consisting of functions which are finely harmonic on the fine interior of $K$. In \cite{DG74} it was shown that these definitions are equivalent. Using a localization result of \cite{BH78} one sees that a function $h\in H(K)$ if and only if it is continuous and finely harmonic on on every fine connected component of the fine interior of $K$. Such collection of sets are usually called {\it restoring}.

Another equivalent definition of $H(K)$ was introduced in \cite{P97} using the notion of Jensen measures which leads another restoring collection of sets.  The main goal of this paper is to reconcile the results in \cite{DG74} and \cite{P97}.

To study these spaces, two notions of Green functions have previously been introduced. One by \cite{P97} as the limit of Green functions on domains $D_j$ where the domains $D_j$ are decreasing to $K$, and alternatively following \cite{F72, F75} one has the fine Green function on the fine interior of $K$. Our Theorem \ref{T:green_equiv} shows that these are equivalent notions.

In Section \ref{S:Jensen} a careful study of the set of Jensen measures on $K$, leads to an interesting extension result (Corollary \ref{C:extend}) for superharmonic functions.  This has a number of applications.  In particular we show that the two restoring coverings are the same.  We are also able to extend some results of \cite{GL83} and \cite{P97} to higher dimensions.
\end{abstract}
\date{\today}
\maketitle
\section{Introduction}

There are several ways to define the spaces ($S(K)$)-$H(K)$ of continuous (super)-harmonic functions  on a compact set $K$ in $\mathbb R^n$. Let $C(K)$ denote the space of all continuous real functions on $K$. The natural definition is to let $H(K)$ or $S(K)$ be the unform closures of  all functions in $C(K)$ which are restrictions of harmonic (resp. superharmonic) functions on a neighborhood of $K$. More fashionably, we can define $H(K)$ and $S(K)$ as the subspaces of $C(K)$ consisting of functions which are finely harmonic (resp. finely superharmonic) on the fine interior of $K$. The equivalence of these definitions was shown in \cite{BH75} and \cite{BH78}.

Another definition was introduced in \cite{P97} using the notion of Jensen measures. A measure $\mu$ supported by $K$ is Jensen  with barycenter $x\in K$ if for every open set $V$ containing $K$ and every superharmonic function $u$ on $V$ we have $u(x)\ge\mu(u)$. The set of such measures will be denoted by $\J_x(K)$. Then $H(K)$ is the subspace of $C(K)$ consisting of functions $h$ such that $h(x)=\mu(h)$ for all $\mu\in\J_x(K)$ and $x\in K$. It was shown in \cite{P97} that this definition is equivalent to the definitions above.

Despite the existence of so many equivalent definitions it is still difficult to verify whether a function on a compact set is harmonic or superharmonic. In \cite{DG74} it was shown that a function $h\in H(K)$ if and only if it is continuous and finely harmonic on the fine interior of $K$. A locatization result from \cite{BH78} implies that a function $h\in H(K)$ if and only if it is continuous and finely harmonic on on every fine connected component of the fine interior of $K$. Such collection of sets are usually called {\it restoring}.

In its turn in \cite{P97} another restoring collection of sets was introduced. For $x\in K$ let $I(x)$ be the set of all points $y\in K$ such that $\mu(V)>0$ for every $\mu\in\J_x(K)$ and every open set $V$ containing $y$. It was shown that the sets $I(x)$ form the restoring covering.

The main goal of this paper is to reconcile the results in \cite{DG74} and \cite{P97}. It required the understanding of a connection between fine topology and Jensen measures. For this we use the fact from \cite{P97} that $I(x)$ is the closure of the set $Q(x)$ of all $y\in K$ such that $G_K(x,y)>0$, where the Green function $G_K$ on $K$ is defined as the limit of Green functions on domains $D_j$ decreasing to $K$.

Fuglede \cite{F72, F75, F83} defined a Green function on $K$ as the fine Green function on the fine interior $int_f(K)$ of $K$. We denote the fine Green function on a finely open set $U$ by
$G^f_U(x,y)$ ( see \cite{F75, F83, F99} for the definition, and
Section \ref{S:basic} for some basic properties).

As the first step we show at Section \ref{S:green} (Theorem \ref{T:green_equiv}) that these two notions of Green functions are constant multiples of each other. These leads to Proposition \ref{P:Pos} which claims that the set $Q(x)$ is a fine connected component of $int_f(K)$.

To finish the reconciliation process in Section \ref{S:Jensen} we study closely the set $\J_x(K)$. The main result (Theorem \ref{T:contraction}) provides Corollary \ref{C:JenRest} claiming that $\mu\in\J_x(K)$ if and only if $\mu\in\J_x(I(x))$.  This corollary proves to be quite useful.  From it we are able to derive a number of applications in Section \ref{S:App}.  In particular Corollary \ref{C:extend} an extension result for superharmonic functions shows that for every $f\in S(I(x))$ there is a $\hat f \in S(K)$ such that $\hat f|_{I(x)} = f$.  Also following from  Corollary \ref{C:JenRest} is the desired reconciliation of the restoring theorem of Poletsky \cite{P97} and the \cite{DG74} result, proved here as Theorem \ref{T:ext_Poletsky}.

In 1983, Gamelin and Lyons have shown \cite[Theorem 3.1]{GL83}
that for $K\subset \mathbb{R}^2$
\[H(K)^\perp = \bigoplus H(\overline{A}_j)^\perp\]
where $A_j$ are the fine components (fine open, fine connected) of the fine interior of $K$.  However their work follows from
an estimate for harmonic measure of the radial projection of a
set, proved by Beurling in his thesis, which has no analog in
$\mathbb{R}^n$ for $n>2$.  By using Theorem \ref{T:ext_Poletsky} we are now able
to extend this result to higher dimensions in Corollary
\ref{C:perp}.  As an application of this we are able to show,
Proposition \ref{P:J=>W}, that every Jensen set is Wermer,
which was first proved by Poletsky in \cite{P97} for $n=2$.

We are especially grateful to Eugene Poletsky for his excellent
guidance and support.

\section{Basic properties}\label{S:basic}
Let $U\subset \mathbb{R}^n$ be a domain.  The set of harmonic
and superharmonic functions on $U$ are denoted by $H(U)$ and
$S(U)$, respectively.  The set of Jensen measures on $U$ at
$x\in U$, denoted $\mathcal{J}_x(U)$, is the collection of all
positive Radon measures $\mu$ with support compactly contained
in $U$ such that $\mu(f)\le f(x)$ for all $f\in S(U)$, where
$\mu(f)=\int fd\mu$. The Jensen measures on $K$ are given by
\[\mathcal{J}_x(K) = \bigcap_{K\subset\, U} \mathcal{J}_x(U), \qquad x\in K\]
where the intersection is over all domains $U$ containing $K$.

Poletsky \cite[Theorem 3.1]{P97} (see also \cite{BH86,BH78,DG74,H85}) has shown a connection
between the set $H(K)$ and the set of
Jensen measures on $K$.
\bT\label{T:Poletsky}
A function is in $H(K)$ if and only if it is continuous and
satisfies the averaging property with respect to every Jensen
measure on $K$, that is
\[H(K)=\left\{h \in C(K)\colon h(x)=\mu(h), \text{ for all } \mu \in \J_x(K) \text{ and every }  x\in K \right\}.\]
\eT
By analogy one can define the set of continuous superharmonic
functions on $K$.  The set of functions $S(K)$ is defined to be the uniform closure of the set of functions $f|_K$ where $f$ is continuous and superharmonic in some neighborhood of $K$.

See \cite{P97} for further properties of the space $H(K)$ and
$\mathcal{J}_x(K)$.  For our purpose it is sufficient to draw
some connections between these sets and their counterparts as
seen from a fine potential theoretic viewpoint.

The two books \cite{B71, F72} are classical references on the
fine topology and many books on potential theory contain
chapters on the topic, e.g. \cite[Chapter 7]{AG01} and
\cite[Chapter 10]{H69}.  Furthermore the topic is generally
subsumed in the more abstract potential theory involving
balayage spaces, e.g. \cite{BH86}.

The \emph{fine topology} on $\mathbb{R}^n$ is the coarsest
topology on $\mathbb{R}^n$ such that all superharmonic functions
are continuous in the extended sense of functions taking values in $[-\infty, \infty]$.  When referring to a topological concept in the fine topology we will follow the standard policy of either using the words ``fine"  or  ``finely" prior to the topological concept or attaching the letter $f$ to the associated symbol. For example, the fine boundary of $K$, $\partial_fK$, is the boundary of $K$ in the fine topology.  The fine topology is strictly finer than the Euclidean topology.

A set $E$ is said to be {\em thin} at a point $x_0$ if $x_0$ is not a fine limit point of $E$, i.e. if there is a fine
neighborhood $U$ of $x_0$ such that $E\setminus\{x_0\}$ does
not intersect $U$.  An example of a thin set is given by the
Lebesgue spine in $\mathbb{R}^3$ defined by
\[L=\{(x, y, z) \colon x>0\, \text{ and }\, y^2+z^2<\exp(-c/x) \},\qquad c>0\]
which is thin at the origin.

Many of the key concepts of classical potential theory have
analogous definitions in relation to the fine topology.
Presently we will recall a few of them.  Relative to a finely
open set $V$ in $\mathbb{R}^n$ the \emph{harmonic measure}
$\delta_x^{\complement V}$ is defined as the swept-out of the
Dirac measure $\delta_x$ on the complement of $V$.  A function
$u$ is said to be \emph{finely hyperharmonic} on a finely open
set $U$ if it is lower finite, finely lower semicontinuous, and
\[ -\infty< \delta_x^{\complement V}(u)\le u(x),\]
 for all $x\in V$ and all relatively compact finely open sets $V$ with fine closure contained in $U$.  A finely hyperharmonic function $u$ is called \emph{finely superharmonic} if $u\not\equiv \infty$, and a function $h$ is said to be \emph{finely harmonic} if $h$ and $-h$ are finely hyperharmonic.  We will need the concept of a fine Green function.  See
\cite{F75, F83, F99} for the definitions and basic properties.  It was shown~\cite{F75} that every
bounded fine open set $U$ admits a fine Green function which we
shall denote by $G^f_U(x,y)$.

Following  \cite{G78} the Choquet boundary of $K$ with respect to $S(K)$ is
\[Ch_{S(K)}K = \{x\in K \colon J_x(K) = \{\delta_x\}\}. \]
Many nice properties of the Choquet boundary are known.  In
particular, we will need the following characterization, see,
for example, \cite[VI.4.1]{BH86} and \cite{H85}.

\bL\label{L:Pfine}  The Choquet boundary of $K$ with respect to $S(K)$
is the fine boundary of $K$, i.e.
\[Ch_{S(K)}K = \partial_f K.\]
\eL
\begin{proof}
Since the fine topology is strictly finer than the Euclidean
topology, any point in the interior of $K$ will also be in the
fine interior of $K$, and any point of $\mathbb{R}^n\setminus
K$  can be separated from $K$ by an Euclidean (therefore fine)
open set.  Therefore the fine boundary of $K$ is contained in
$\partial K$. The result follows immediately from \cite[Theorem
3.3]{P97} or \cite[Proposition 3.1]{BH86} which states that
$\mathcal{J}_x(K) = \{\delta_x\}$ if and only if the complement
of $K$ is non-thin at $x$, i.e. $x$ is a fine boundary point of
$K$.
\end{proof}
\par In particular,
\bC\label{C:loc}
If $\J_x(K) \neq \{\delta_x\}$, then $x \in int_fK$.
\eC

The set $\partial_fK$ is also called the stable boundary of
$K$.  In fact the lemma shows that $Ch_{S(K)}K$ is the finely
regular boundary of the fine interior of $K$.  For more details
on finely regular boundary points and other related concepts,
see \cite[VII.5-7]{BH86} and \cite{H85}.

The following lemma has been known since the book of
Fuglede~\cite[p. 147]{F72}.
\bL
A fine open set $U$ in $\mathbb{R}^n$ has at most countably
many fine open connected components.
\eL

\section{On the Green function associated to a compact set}\label{S:green}
We now proceed to study the Green function on $K$.
Recall~\cite[Theorem (b) 1.VII.6, p. 94]{D83} that if $D$ is an
open Greenian set in $\mathbb{R}^n$ so that $\{D_j\}$ is a
decreasing sequence of open sets converging to $D$, then the
sequence $\{G_{D_j}(\cdot,y)\}$ of Green functions associated
to $\{D_j\}$ is decreasing to $G_D(\cdot,y)$ for every $y\in
D$.  By analogy one can define a Green function on a compact
set $K$ as the limit of the sequence $\{G_{D_j}(\cdot,y)\}$
where $y\in K$ and $\{D_j\}$ is any decreasing sequence of open
sets converging to $K$.  In the article~\cite{P97} Poletsky
defines a Green function on a compact set in this way.

Recall, \cite[p. 90]{D83}, that for a regular open set $D$ the
associated Green function $G_D(\cdot,y)$ extends continuously
as $\hat{G}_D(\cdot,y)$ to $\mathbb{R}^n$ for any $y\in D$
where $\hat{G}_D(\cdot,y)=0$ on $\complement D$, the complement
of $D$, and this extension $\hat{G}_D(\cdot,y)$ is subharmonic
on $\mathbb{R}^n\setminus\{y\}$.

In the following proposition we outline some of the basic
properties of $\hat{G}_K$.

\bP\label{P:pgk} For all $y\in K$, the function
$\hat G_K(\cdot,y)\colon \mathbb{R}^n\rightarrow [0, \infty]$
defined as $\hat G_K(\cdot,y) = \liminf_j
\hat{G}_{D_j}(\cdot,y)$ has the following properties:
\begin{enumerate}
\item[{\rm i.}]\label{L:ExtBy0}  $\hat G_K(x, y)=0$ when
    $x\in \complement K := \mathbb{R}^n\setminus K$ and
    $y\in K$,
\item[{\rm ii.}] $\hat G_K$ does not depend on the sequence
    $\{D_j\}$ chosen,
\item[{\rm iii.}] $\hat G_K \ge 0$ and $\hat G_K(y,
    y)=+\infty$ for all $y\in K$,
\item[{\rm iv.}] $\hat G_K$ is symmetric, i.e. $\hat G_K(x,
    y) = \hat G_K(y, x)$, for all $x,y\in K$,
\item[{\rm v.}] $\hat G_K(\cdot,y)$ is super-averaging on
    $K$, i.e. $\hat G_K(x, y) \ge \int \hat G_K(\zeta, y)\,
    d\mu(\zeta)$ for all $\mu\in \J_x(K)$ with with $x \in
    K$, and
\item[{\rm vi.}]\label{L:sub} $\hat G_K(\cdot,y)$ is
    subharmonic on $\mathbb{R}^n\setminus\{y\}$.
\end{enumerate}
\eP

\begin{proof}[proof of i.]
This follows from the fact that $\hat G_{D_j}(x,y)=0$ whenever
$x\notin \overline{D}_j$.
\end{proof}
\begin{proof}[proof of ii.]
If $D_1\supset D_{2}$ then $\hat{G}_{D_1}(\cdot,y)\ge
\hat{G}_{D_{2}}(\cdot,y)$ for any Greenian sets $D_1$ and $D_2$.
Alternatively we could have defined $\hat{G}_K$ by \[\hat{G}_K(\cdot,y) =
\inf\{ \hat{G}_D(\cdot,y)\colon D\supset K, D \text{
Greenian}\},\qquad y\in K.\qedhere\]
\end{proof}
\begin{proof}[proof of iii.] As $\hat{G}_D \ge 0$ and $\hat{G}_D(y, y)=+\infty$ for all
$x,y \in D$ for any Greenian $D$.
\end{proof}
\begin{proof}[proof of iv.] Since $\hat{G}_D(x, y) = \hat{G}_D(y, x)$ for all
 $x,y \in D$ for any Greenian $D$.
\end{proof}
\begin{proof}[proof of v.] For any Greenian set $D$ the function
$\hat{G}_D(\cdot,y)$ is superharmonic on $D$.  Then $\hat{G}_D(x, y) \ge
\int \hat{G}_D(\zeta, y)\, d\mu(\zeta)$ for all $\mu\in \J_x(D)$ with
$\zeta \in D$.  If $D_j$ is a decreasing sequence of domains
converging to $K$, then $\hat{G}_{D_j}(\cdot,y)$ is decreasing to
$\hat{G}_K(\cdot,y)$.  Therefore by the Lebsegue Monotone Convergence
Theorem $\hat{G}_K(x,y) \ge \int \hat{G}_K(\zeta, y)\, d\mu(\zeta)$ for all
$\mu\in \cap_j \J_x(D_j):=\J_x(K)$ with $x \in K$.
\end{proof}
\begin{proof}[proof of vi.] Let $\{D_j\}$ be a decreasing sequence of
regular domains converging to $K$.  Then
$\hat{G}_{D_j}(\cdot,y)$ is continuous, and so $\hat{G}_K(\cdot,y)$
must be upper semicontinuous.  For any $j$ and any $y\in D_j$,
by \cite[p. 90]{D83} the extension $\hat{G}_{D_j}(\cdot,y)$ of
Green function $\hat{G}_{D_j}(\cdot,y)$ by $0$ is subharmonic on
$\mathbb{R}^n\setminus \{y\}$.  Therefore by the Lebesgue Monotone
Convergence Theorem $\hat{G}_K(\cdot,y)$ is subaveraging on
$\mathbb{R}^n\setminus \{y\}$ as it is the decreasing limit of
a sequence of subharmonic functions.  Since $\hat{G}_K(\cdot,y)$ is
upper semicontinuous and subaveraging, $\hat{G}_K(\cdot,y)$ is
subharmonic on $\mathbb{R}^n\setminus \{y\}$.
\end{proof}

It was shown in \cite{F75} that every bounded fine open set
$U$ admits a fine Green function which we shall denote by
$G^f_U(x,y)$.  The following result shows that for a compact
set $K$ the functions $\hat{G}_K(x,y)$ and $G^f_{int_fK}(x,y)$ are
scalar multiples of each other.

\bT\label{T:green_equiv}
For any compact set $K\subset \mathbb{R}^n$, $n\ge 2$, there is
$c>0$ such that $\hat{G}_K(x,y)=cG^f_{int_fK}(x,y)$ for any $y\in
int_fK$.
\eT
\begin{proof}
Fuglede has given a simple characterization of the fine Green
function up to multiplication by a positive constant.  Indeed,
if a function $g\colon U\times U\rightarrow \mathbb{R}$ has the
following properties
\begin{enumerate}
\item[{\rm 1.}] $g(\cdot,y)$ is a nonnegative finely
    superharmonic function on $U$,
\item[{\rm 2.}] if $v$ is finely subharmonic on $U$ and
    $v\le g(\cdot,y)$, then $v\le 0$,
\item[{\rm 3.}] $g(\cdot, y)$ is finely harmonic on
    $U\setminus\{y\}$ for any $y\in U$, and
\item[{\rm 4.}] $g(y,y)=+\infty$
\end{enumerate}
then $g(x,y)=cG^f_U(x,y)$ for some $c>0$ for all $x,y\in U$.

Hence to prove the theorem we need only to check these
properties. Firstly, we note that by Lemma \ref{L:sub}
$\hat{G}_K(\cdot,y)$ is subharmonic (and thereby finely subharmonic)
on $\mathbb{R}^n\setminus\{y\}$, which implies (\cite[Theorem
9.10]{F72}) fine continuity on $\mathbb{R}^n\setminus\{y\}$.

In fact, we will shall now see that $\hat{G}_K(\cdot,y)$ is finely
continuous at $y$ when $y\in int_fK$. Every bounded fine open
set admits a fine Green function, cf.~\cite{F75, F99}.  Let
$G^f_{int_f K}$ denote the fine Green function corresponding to
the bounded fine open set $int_fK$. Since $int_f K\subset D_j$
we have $G^f_{int_f K}(\cdot,y) \le \hat{G}_{D_j}(\cdot,y)$.  As
$\hat{G}_K(\cdot,y)$ is the decreasing limit of $\hat{G}_{D_j}(\cdot,y)$ we
have the inequalities
\[G^f_{int_f K}(\cdot,y) \le \hat{G}_K(\cdot,y) \le \hat{G}_{D_j}(\cdot,y),\]
for all $y\in int_fK$.   Since $G^f_{int_f K}(\cdot,y)$ and
$\hat{G}_{D_j}(\cdot,y)$ are finely continuous, $\hat{G}_K(\cdot,y)$ must
be finely continuous at $y$ as
\[\infty =  \text{f -}\lim_{x\rightarrow y}G^f_{int_f K}(x, y)\le \text{f -}\lim_{x\rightarrow y}\hat{G}_K(x, y)\le \text{f -}\lim_{x\rightarrow y}\hat{G}_{D_j}(x, y)=\infty.\]
Therefore $\hat{G}_K(\cdot,y)$ is finely continuous on $\mathbb{R}^n$
when $y\in int_fK$.

Thus $\hat{G}_K(\cdot,y)$ is finely superharmonic on $int_f K$ as it
is finely continuous and the decreasing limit of
$\{\hat{G}_{D_j}(\cdot,y)\}$, a sequence of finely
superharmonic functions on $int_f K$ and this implies that 1.
holds.

Suppose that $\hat{G}(x_0,y)>0$ for $x_0\in \partial_fK$. Then there
is a fine neighborhood $V$ of $x_0$ such that $\hat{G}_K(x,y)>0$ for
all $x\in V$. By definition $x\in \partial_f K$ if and only if
$\complement K$ is non-thin at $x$. As $x_0\in \partial_fK$,
this means that $V\cap \complement K \neq \emptyset$.  However
by Lemma \ref{L:ExtBy0}, $\hat{G}_K(x, y)=0$ for $x\in \complement K$
and $y\in K$, a contradiction.  Therefore $\hat{G}_K(x,y)=0$ for all
$x\in\partial_fK$ and $y\in int_fK$. So $\hat{G}_K$ is a fine
potential on $int_f K$ by the minimum principle
\cite[III.4.1]{BH86} (see also \cite[Theorem 9.1]{F72}) and
this implies 2.

We have seen above that $\hat{G}_K(\cdot,y)$ is finely superharmonic
on $int_fK$. By Proposition  \ref{P:pgk}.vi $\hat{G}_K(\cdot,y)$ is
finely subharmonic on $int_fK\setminus\{y\}$. Therefore
$\hat{G}_K(\cdot,y)$ is finely harmonic on $int_fK\setminus\{y\}$ and
we checked 3.

The property 4. follows immediately from Proposition
\ref{P:pgk}.iii and the theorem is proved.
\end{proof}

\bP\label{P:Pos}
The Green function $\hat{G}_K(x, y) > 0$ for $x, y \in K$ if and only
if $x$ and $y$ are in the same fine connected component of
$int_fK$.
\eP
\begin{proof}
By the previous proposition $\hat{G}_K(\cdot,y)$ is finely
superharmonic on $int_f K$.  If $\hat{G}_K(x, y)=0$, then by
\cite[Theorem 12.6]{F72} for all $\zeta$ in the fine component
of $y$ we have $\hat{G}_K(\zeta, y)=0$.  Therefore $\hat{G}_K(\cdot,y)>0$
on the fine component containing $y$.

Suppose that $int_f K$ has multiple components.  Each component
is fine open and therefore has its own Green function.  We can
define a function $g(x,y)$ on $int_f K$ by
\[ g(x, y)= \begin{cases} G^f_{Q_x}(x, y), & y\in Q_x\\ 0, & y\in (int_f K)\setminus Q_x\end{cases}\]
where $Q_x$ is the fine component containing $x$.  Since fine
subharmonicity and fine harmonicity are local properties, $g$
satisfies the requirements mentioned in the proof of Theorem \ref{T:green_equiv} to be a
positive multiple of the fine Green function on $int_fK$.
Therefore $G_{int_fK}^f(x,y)$ is positive if and only if $x$
and $y$ are in the same fine component of $int_fK$.  So $\hat{G}_K(x,
y)=0$ when $x$ and $y$ are in different fine connected components.

In the proof of the previous proposition we proved that $\hat{G}_K(x,
y)=0$ for $x\in \partial_f K$ and $y\in K\setminus \{x\}$.
\end{proof}

In \cite{P97} Poletsky introduced the sets
\[Q(x)=\{y\in K\colon \hat{G}_K(x, y)>0\},\]for every $x\in K$.  The following corollary directly
follows from Proposition \ref{P:Pos} and characterizes these sets
in terms of the fine topology.

\bC\label{C:f-comp}
For all $x\in int_fK$, the set $Q(x)$ is the
fine connected component of $int_fK$ which
contains $x$. Additionally the point $x\in K$ is in $\partial_f
K$ if and only if $Q(x)=\{x\}$.
\eC

\section{Jensen measures}\label{S:Jensen}

Some results from~\cite{P97} now follow from standard
properties of the fine potential theory and the fine topology.
For example \cite[Theorem 3.6 (2)]{P97} is the partitioning the set
$K$ into the fine connected components of $int_fK$ and
singleton sets for peak points (i.e. the set $\partial_fK$) forms an equivalence relation,
\cite[Theorem 3.6 (3)]{P97} is the fine minimum principle, and
\cite[Theorem 3.6 (4)]{P97} is that fine connected components have
positive measure.  We can now extend/rephrase some results
of~\cite{P97} and use them to obtain some new results.

\bT\label{T:PJenVis}
For $x\in K$ and any $\varepsilon>0$ there exists a $\mu\in
\J_x(K)$ with $\mu(B(y, \varepsilon))>0$ if and only if the
point $y$ is in the (Euclidean) closure $\overline{Q(x)}$ of
the fine component of $x$.
\eT
\begin{proof}
In \cite{P97} Poletsky defines $I(x)$ as the set of points
$y\in K$ with the property that for any $\varepsilon>0$ there
exists a $\mu\in \J_x(K)$ with $\mu(B(y, \varepsilon))>0$ and
in \cite[Theorem 3.6 (1)]{P97} proves that $I(x) =
\overline{Q(x)}$.  The result follows from Corollary
\ref{C:f-comp}.
\end{proof}

The following corollary is an immediate consequence of the
previous theorem.
\bC\label{C:JenSupp}
Let $K$ be a compact set in $\mathbb{R}^n$, $n\ge 2$.  Then
$\supp \mu \subset \overline{Q(x)}$ for all $\mu \in \J_x(K)$.
\eC

For use in the following proposition we recall the notion of a
reduced function, see \cite[Definition 5.3.1]{AG01}.  Fix a
Greenian open set $\Omega\subset\mathbb{R}^n$. Let
$U_+(\Omega)$ be the set of non-negative superharmonic
functions on $\Omega$.  For $u\in U_+(\Omega)$ and $E\subset
\Omega$, the reduced function of $u$ relative to $E$ in
$\Omega$ is defined by
\[R_u^E(x) = \inf\{v(x)\colon v\in U_+(\Omega) \text{ and }
v\ge u \text{ on } E\},\qquad x\in \Omega.\]
Also note that $\hat{R}_u^E$ is the lower semicontinuous
regularization of $R_u^E$.

\bP\label{P:polar}
Let $U$ and $V$ be disjoint fine open sets.  Then
$V\cap\overline{U}$ is a polar set.
\eP
\begin{proof} It suffices to prove this statement when $U$ and
$V$ are bounded. Otherwise, we may consider intersections of
these sets with increasing sequence of open balls.

Let $\Omega$ be any open Greenian set containing $U$ and $V$.
Since $U$ is disjoint from $V$, $U$ is thin at $y$ for every
$y\in V$. Then by \cite[Theorem 7.3.5]{AG01} there is a bounded
continuous potential $u^\#$ on $\Omega$ with the property that
$\hat{R}^{U}_{u^\#}(y)<u^\#(y)$ for all $y\in V\cap
\overline{U}$. By construction $R^{U}_{u^\#} \ge u^\#$ and
$R^{U}_{u^\#}(x) = \hat{R}^{U}_{u^\#}(x)=u^\#(x)$ for all $x\in
U$.  Therefore $V\cap\overline{U} \subset \{R^{U}_{u^\#} \neq
\hat{R}^{U}_{u^\#}\}$, and by \cite[Theorem 5.7.1]{AG01} the
set $\{R^{U}_{u^\#} \neq \hat{R}^{U}_{u^\#}\}$ is polar.
\end{proof}

\bC\label{C:f-int}
For a compact set $K\subset \mathbb{R}^n$, let $\{A_i\}$ be the
collection of disjoint fine connected components of the fine interior of
$K$.  Then $int_f\, \overline{A_i}=A_i$ for all $i$.
\eC
\begin{proof}
We will show that $int_f\, \overline{A_i}$  has only one fine
component and so it must be $A_i$.  Suppose that
$int_f\,\overline{A_i} = A \cup V$ where $A$ is the fine
component containing $A_i$ and $V$ is fine open and disjoint
from $A$.  First we note that $A_i=A$ as $A_i\subset A\subset
int_fK$ and $A_i$ is a fine component of $int_fK$.  Secondly,
$V$ is disjoint from $A_i$ and contained in
$int_f\,\overline{A_i}$, hence $V\subset \overline{A_i}\setminus
A_i$.  Therefore by Proposition \ref{P:polar}, we have that $V$
must be polar and cannot be fine open.
\end{proof}

The following corollary tells us that the only the trivial
Jensen measures can have support in the closure of two fine
components.  We use the notation $\J(K):=\cup_{x\in K}\J_x(K)$ to
denote the collection of all Jensen measures on $K$.
\bC
Let $\{A_j\}$ be the fine connected components of the fine interior of
$K$.  Then
\[\J(\overline{A}_i) \bigcap \J(\overline{A}_j) = \bigcup_{x\in \overline{A}_i\cap \overline{A}_j}\{\delta_x\},\]
where $i\neq j$.
\eC
\begin{proof}
Let $\mu \in \J(\overline{A}_i) \bigcap \J(\overline{A}_j)$
with $i\neq j$.  Then there is an $x_i \in \overline{A}_i$ and
$x_j \in \overline{A}_j$ so that $\mu \in
\J_{x_i}(\overline{A}_i) \bigcap \J_{x_j}(\overline{A}_j)$.  As
the coordinate functions are harmonic, this implies that
$x_i=x_j$.  Let us call $x_0:=x_i=x_j \in \overline{A}_i
\bigcap \overline{A}_j$.  As $A_i$ and $A_j$ are disjoint, we
have by Corollary \ref{C:f-int} that $x_0$ must be in the fine
boundary of either $\overline{A}_i$ or $\overline{A}_j$.
However the only way that $x_0$ can be in the fine boundary
(see Lemma \ref{L:Pfine}) is if $\mu=\delta_{x_0}$.
\end{proof}

The following theorem gives sufficient condition on a subset $E$ of $K$ so that the Jensen measures on $K$ with barycenter $x\in E$ belong to the Jensen measures on $E$.

\bT\label{T:contraction}
Let $A\subset K \subset \mathbb{R}^n$, $n\ge 2$, with $K$
compact with $A$ and $int_fK\setminus A$ fine open, that is $A$ is the union of fine connected components of $int_fK$.  Suppose that $\supp \mu \subset \overline{A}$ for all $\mu\in
\J_x(K)$ and all $x\in A$ then $\J_x(K)\subset
\J_x(\overline{A})$ for all $x\in A$.
\eT
\begin{proof}
Suppose there exists $\mu\in \J_{x_0}(K)$, $f\in
S(\overline{A})$ and $a>0$, such that $\mu(f)>f({x_0})+a$.  As
$cf+c'$ is also in $S(\overline{A})$ for $c>0$ and since the
functions in $S(\overline{A})$ are uniform limits of continuous
superharmonic functions defined in neighborhoods of
$\overline{A}$, we may assume that $f\in C(G)\cap S(G)$ for
some open set $G\supset \overline{A}$ with the properties
$-1<f({x_0})<\mu(f)<0$ and $-1<f<0$.  Let
$a:=\mu(f)-f({x_0})>0$ and take $G'$ open with
$\overline{A}\subset G'$ and $\overline{G'}\subset G$.

Pick $\phi\in C(\mathbb{R}^n)$ with $\phi =0$ on
$\overline{A}$, $\phi=1$ on $\mathbb{R}^n\setminus G'$ and $0<
\phi < 1$ on $G'\setminus \overline{A}$.  By Edwards Theorem
(see \cite{CR97})
\[E\phi(y)=\inf\{f(y)\colon f\in S(K), \; f\ge \phi\} = \sup\{\nu(\phi)\colon \nu\in \J_y(K)\}.\]
By assumption $\supp(\nu)\subset \overline{A}$ for all
$\nu \in \J_y(K)$ and every $y\in A$. So $E\phi(y)=0$ for
every $y\in A$. Therefore there exists a $g\in S(K)$ with $0\le
\phi\le g \le 1$ and $g(x_0)<\varepsilon < a/3$.

Actually we can say a little more.  By Corollary \ref{C:loc},
we know that $\J_y(K)\neq \{\delta_y\}$ if and only if $y\in
int_fK$. This allows us to decompose $\overline{A}$ into three
sets; $A$, $\partial_1 A\subset \partial A$ where $\J_y(K)
=\{\delta_y\}$ for $y\in \partial_1A$, and $\partial_2A =
\overline{A} \setminus (A \cup \partial_1 A)$.  Each point in
$\partial_2A$ belongs to $int_fK \setminus A$.  Recall that by hypothesis $int_fK\setminus A$ is fine open.  Therefore $\partial_2A \subset \overline{A}\cap (int_fK\setminus A)$, which means that $\partial_2A$ is polar by Proposition \ref{P:polar}. Since $\partial_2A$ is a polar set, we see that $\mu(\partial_2A)=0$.

Thus there exists $C$ a compact neighborhood of $x$ with
$C\subset A\cup \partial_1A$ so that $\mu(C)>1-\varepsilon$.
As $\left.E\phi\right|_{A\cup \partial_1A}=0$, trivially
$\left.E\phi\right|_C=0$.  For every $y\in C$ there are a
continuous and superharmonic function $g_y\ge\phi$ in a neighborhood
of $K$ and  an open neighborhood $U_y$ of $y$ with $g_y<\varepsilon$ on $U_y$.
The sets $U_y$ cover $C$, so by compactness we can pick up $y_1, \ldots
, y_N$ so that $C\subset U_{y_1} \cup \cdots \cup U_{y_N}$.
Then $g = \min\{g_{y_1}, \ldots , g_{y_n}\}$ has the property
$g|_C<\varepsilon$  and $\mu(\{g > \varepsilon\})<\varepsilon$.

Consider the function $f+g$.  As $g\ge 0$ we have
\[\mu(f+g) =
\mu(f)+\mu(g) >f(x_0) +g(x_0) +a-g(x_0) > (f(x_0) + g(x_0))
+a-\varepsilon.\]  As $\phi\le g$, we have that $f+g\ge 0$ on
$K\setminus \overline{G'}$.  Note also that
\[f(x_0)+g(x_0)=\mu(f)-a+g(x_0)<-a+g(x_0)<-a+\varepsilon<0.\] So
\[h(y) = \begin{cases}0, & K\setminus \overline{G}\\
\min\{f+g,0\}, & G\cap K\end{cases}\] is in $C(K)$, $h\equiv0$ on $K\setminus G'$ and $h(x_0) = f(x_0)+g(x_0)$.

To see that $h$ is in $S(K)$ we use a localization argument.
Let $\mathcal{V}$ be a covering of the fine interior of $K$ by
fine open sets such that $V\in \mathcal{V}$ has the property:
if $V\cap G' \neq \emptyset$ then $V\subset G$.  If $V\subset
G$, then $h=\min\{f+g, 0\} \subset S(V)$.  If $V\cap G'=\emptyset$
then $h\equiv 0 \in S(V)$.  Thus
$h\in S(K, int_fK, \mathcal{V} )= S(K)$, by \cite[Proposition
3.5]{BH78}.

Thus \[\mu(h)=\int\limits_{\{f+g<0\}}(f+g)\,d\mu = \mu(f+g) -
\int\limits_{\{f+g\ge0\}}(f+g)\,d\mu.\] Now
$\mu(f+g)>f(x_0)+g(x_0)+a-\varepsilon$ and
\[\int\limits_{\{f+g\ge0\}}(f+g)\,d\mu =\int\limits_{\{f+g\ge0\}}f\,d\mu
+\int\limits_{\{f+g\ge0\}}g\,d\mu.\] The first integral on the
right is negative and because $g\ge0$
\[\int\limits_{\{f+g\ge0\}}g\,d\mu\le \int g\,d\mu.\]
But the last integral is equal to
\[\int\limits_{\{g>\varepsilon\}}g\,d\mu+
\int\limits_{\{g\le\varepsilon\}}g\,d\mu\le\mu(\{g>\eps\})+\varepsilon.\]
Recall that $\mu(\{g>\eps\})<\varepsilon$. Thus
\[\mu(h)>f(x_0)+g(x_0)+a-3\varepsilon>h(x_0).\]
However this contradicts that $\mu\in J_{x_0}(K)$ and $h\in
S(K)$. Hence $\J_x(K)\subset\J_x(\overline{A})$ for all $x\in
A$. Since the reverse inclusion is trivial the theorem is
proved.

\end{proof}

We also get the following useful restriction property of Jensen measures.
\bC\label{C:JenRest}
Let $K\subset\mathbb{R}^n$, $n\ge 2$, be a compact set.  For
all $x\in K$, we have
\[\J_x(\overline{Q(x)})=\J_x(K).\]
\eC
\begin{proof}
We will show $\J_x(\overline{Q(x)})\subset \J_x(K)$ first.
Consider $\mu\in \J_x(\overline{Q(x)})$ and $u\in S(K)$.  Then
$u|_{\overline{Q(x)}}\in S(\overline{Q(x)})$, so that $u(x)\ge
\mu(u)$.  Thus $\mu\in \J_x(K)$.

By Corollary \ref{C:JenSupp}, $\supp \mu \subset
\overline{Q(x)}$ for all $\mu\in \J_x(K)$, which by Theorem
\ref{T:contraction} means that $\J_x(K)\subset
\J_x(\overline{Q(x)})$.
\end{proof}

\section{Applications}\label{S:App}

An interesting corollary follows immediately
from the proof of Theorem \ref{T:contraction}.
For any cone of functions $\mathcal{R}$, we define the closure $\widetilde{\mathcal{R}}$ of $\mathcal{R}$ as all continuous functions which can be represented as the infimum of functions from $\mathcal{R}$.
\bC\label{C:extend}
Let $K\subset\mathbb{R}^n$, $n\ge 2$, be a compact set with
$\{A_j\}$ the fine connected components of the fine interior of $K$.
Then
\[S(\overline{A}_j) =\widetilde{\left.S(K)\right|_{\overline{A}_j}}\]
for every component $A_j$.
\eC
\begin{proof}
It is clear that $\left.S(K)\right|_{\overline{A}_j}\subset S(\overline{A}_j)$.  Consider any function $f\in S(\overline{A}_j)$.
By Edwards Theorem (see \cite{CR97})
\[ f(x) = \inf\{\phi(x)  \colon\phi \in S(\overline{A}_j) \text{ and } \phi \ge f \text{ on } \overline{A}_j\} = \sup\{\mu(f) \colon \mu \in \J_x(\overline{A}_j)\},\]
for all $x\in \overline{A}_j$.  From Corollary \ref{C:JenRest} we have $\J_x(\overline{A}_j) = \J_x(K)$ when $x\in A_j$.  Therefore we may apply Edwards Theorem again to see that
\[ f(x) = \sup\{\mu(f) \colon \mu \in \J_x(K)\}= \inf\{\phi(x) \colon\phi \in S(K) \text{ and }  \phi \ge f \text{ on }  K\},\]
for $x\in A_j$.  Thus $f\in \widetilde{\left.S(K)\right|_{\overline{A}_j}}$.
\end{proof}

The following theorem shows that the restoring covering of \cite{P97} is given by the fine connected components of $int_fK$.

\bT\label{T:ext_Poletsky}
Let $K\subset \mathbb{R}^n$, $n\ge2$, be a compact set with
$\{A_j\}$ denoting the fine components (fine open, fine
connected) of the fine interior of $K$.  For any $f\in C(K)$,
$f \in H(K)$ if and only if $f \in H(\overline{A}_j)$ for all
$j$.
\eT
\begin{proof}
Recall that
\begin{align*}
H(K) = \{f \in C(K) \colon f(x) = \mu(f) \text{ for all } \mu \in \J_x(K) \text{ and every } x\in K\}
\end{align*}
However if $x$ is in $A_j$ by Corollary \ref{C:JenRest} we have that $\J_x(\overline{A}_j)=\J_x(K)$ which implies the result.
\end{proof}

As a corollary we may extend the \cite{GL83} result to higher
dimensions.  Recall that for any compact set $E$, the set
$H(E)^\perp$ is the set of Radon measures $\mu$ with
$supp(\mu)\subset E$ such that $\mu(h)=0$ for all $h\in H(E)$,
and if $m(E)$ is any set of Radon measures with support in $E$
the set $^\perp m(E)$ consists of all $f\in C(E)$ such that
$\mu(f)=0$ for all $\mu\in m(E)$.

\bC\label{C:perp}
For any $K\subset \mathbb{R}^n$, $n\ge 2$, compact
\[H(K)^\perp = \bigoplus H(\overline{A}_j)^\perp\]
where $A_j$ are the fine components (fine open, fine connected)
of the fine interior of $K$.
\eC
\begin{proof}
Consider any $\mu \in \oplus H(\overline{A}_j)^\perp$ and $h\in
H(K)$.  Then $h|_{\overline{A}_j} \in H(\overline{A}_j)$, so
$\mu(h)=0$.  Thus $\oplus H(\overline{A}_j)^\perp \subset
H(K)^\perp$.

Conversely, suppose that $h\in C(K)$ and $h\in \
^\perp\!\left(\oplus H(\overline{A}_j)^\perp\right)$.  Then
$h|_{\overline{A}_j} \in \
^\perp\!\left(H(\overline{A}_j)^\perp\right) =
H(\overline{A}_j)$.  The restoring property (Theorem \ref{T:ext_Poletsky}) then
implies that $h\in H(K)$.  Therefore $ \ ^\perp\!\left(\oplus
H(\overline{A}_j)^\perp\right) \subset H(K)$ and so $H(K)^\perp
\subset \oplus H(\overline{A}_j)^\perp$.
\end{proof}

Recall the following definitions of Poletsky \cite[Def 3.9,
3.15]{P97}.
\bD A compact set $K\subset \mathbb{R}^n$, $n\ge 2$, is called Jensen if
$K=\overline{Q(x)}$ for some $x\in K$, and Wermer if for all
$x\in K$, either $\overline{Q(x)}=K$ or
$\overline{Q(x)}=\{x\}$.
\eD

It has been shown in \cite[Corollary 3.16]{P97} that every
Jensen set is a Wermer set in the plane.  We can now provide a
proof of this in $\mathbb{R}^n$.

\bP\label{P:J=>W}
A Jensen set is Wermer.
\eP
\begin{proof}
Suppose $K$ is Jensen.  Then $K = \overline{Q(x_0)}$ for some
$x_0\in K$.  Every $y\in K$ is either a fine boundary point or
in the fine interior.  If $y$ is in the fine boundary of $K$,
then $Q(y) = \{y\}$ by Corollary \ref{C:f-comp}.

We will show that $int_fK$  has only one fine component which
must be $Q(x_0)$.  Suppose that $int_fK = Q(x_0) \cup V$ were
$V$ is fine open and disjoint from $Q(x_0)$.  Since
$\overline{Q(x_0)}=K$, by Proposition \ref{P:polar}, we have
that $V$ must be polar and cannot be fine open.

Thus for any $y \in int_fK$, we have $Q(y)=Q(x_0)$ and so
$\overline{Q(y)}=K$.
\end{proof}

The set $K = [0, 1] \subset \mathbb{R}^2$ provides a simple
example of a Wermer set that is not Jensen.  Every point is a
fine boundary point, so $Q(x) = \{x\}$ for all $x\in K$.
However there is no point $x_0\in K$ such that
$K=\overline{Q(x_0)}$.  Proposition \ref{P:J=>W} can be
interpreted as saying that if $K$ is Wermer then either
$H(K)=C(K)$ or $K$ is Jensen.


\begin{thebibliography}{XXXX}
\bibitem[AG01]{AG01} Armitage, D. H., Gardiner, S. J., {\em
    Classical Potential Theory,} Springer-Verlag, London, 2001.
\bibitem[BH86]{BH86}  Bliedtner, J., Hansen, W., {\em Potential
    Theory,} Springer-Verlag, 1986.
\bibitem[BH75]{BH75} Bliedtner, J., Hansen, W., {\em Simplicial
    Cones in Potential Theory,} Inventiones Mathematicae,
    Volume 29, (1975), 83--110
\bibitem[BH78]{BH78} Bliedtner, J., Hansen, W., {\em Simplicial
    Cones in Potential Theory II (Approximation Theorems),}
    Inventiones Mathematicae, Volume 46, (1978), 255--275
\bibitem[B71]{B71} Brelot, M., {\em On topologies and
    Boundaries
    in Potential Theory,} Lecture Notes in Mathematics 175,
    Springer-Verlag, Berlin, 1971.
\bibitem[CR97]{CR97} Cole, B. J., Ransford, T. J., {\em
    Subharmonicity
    without Upper Semicontinuity,} Journal of Functional
    Analysis, Volume 147, (1997), 420--442
\bibitem[DG74]{DG74} Debiard, A., Gaveau, B., {\em Potential
    fin et
    alg\'{e}bres de fonctions analytiques,} I.J. Funct.  Anal.,
    Volume 16, (1974), 289--304
\bibitem[D83]{D83}  Doob, J. L., {\em Classical Potential
    Theory and Its Probabilistic Counterpart,}  Grundlehren der
    mathematischen Wissenschaften 262, Springer-Verlag, 1983.
\bibitem[F72]{F72} Fuglede, B., {\em Finely Harmonic
    Functions,}
    Lecture Notes in Mathematics 289, Springer-Verlag, 1972.
\bibitem[F75]{F75} Fuglede, B., {\em Sur la fonction de Green
    pour un domaine fin,} Annales de l'institiut Fourier, tome
    25, no 3-4 (1975), 201--206.
\bibitem[F83]{F83} Fuglede, B., {\em Integral Representation of
    Fine Potentials,} Math. Ann., 262 (1983), 191--214.
\bibitem[F99]{F99} Fuglede, B., {\em The Dirichlet Laplacian on
    finely open sets.} Potential Anal., 10 (1999), no. 1,
    91--101.
\bibitem[G78]{G78} Gamelin, T. W., {\em Uniform
    Algebras and Jensen Measures,} Cambridge University Press, 1978
\bibitem[GL83]{GL83}  Gamelin, T. W., Lyons, T. J., {\em Jensen
    Measures for R(K),} J. London Math. Soc., Volume 27, Number
    2,  (1983), 317--330.
\bibitem[H85]{H85} Hansen, W., {\em Harmonic and superharmonic
    functions on compact sets.} Illinois J. Math., Volume 29,
    Number 1, (1985), 103--107.
\bibitem[H69]{H69} Helms, L. L., {\em Introduction to
    Potential Theory,} Wiley-Interscience, 1969.
\bibitem[H62]{H62} Herv\'{e}, R.-M., {\em Recherches
    axiomatiques sur la th\'{e}orie des fonctions
    surharmoniques et du potentiel,} Ann. Inst. Fourier, 12,
    (1962), 415--517.
\bibitem[P97]{P97} Poletsky, E. A., {\em Approximation by
    Harmonic Functions,} Transactions of the AMS, Volume 349,
    Number 11, (1997), 4415--4427.
\end{thebibliography}
\end{document}